\documentstyle[amstex,amssymb]{article}

\renewcommand{\phi}{\varphi}

\newtheorem{theorem}{Theorem}[section]
\newtheorem{proposition}[theorem]{Proposition}
\newtheorem{lemma}[theorem]{Lemma}
\newtheorem{corollary}[theorem]{Corollary}
\newtheorem{remark}[theorem]{Remark}
\newtheorem{definition}[theorem]{Definition}
\newtheorem{question}[theorem]{Question}
\newtheorem{example}[theorem]{Example}

\begin{document}
\thispagestyle{empty}

\begin{center}
{\Large\bf Dual Entropy in Discrete Groups with Amenable Actions}
\end{center}

\begin{center}{\bf Nathanial P. Brown \footnote{ 
      NSF Postdoctoral Fellow. 

      1991 AMS Classification; 46L55.}}\\
   UC-Berkeley\\
   Berkeley, California 94720 \\
   nbrown{\char'100}math.berkeley.edu 

{\bf Emmanuel Germain}

Universit\'e de Paris VII \\
UFR de Math\'ematiques, 75005 Paris \\
germain{\char'100}math.jussieu.fr

\end{center}

\begin{abstract}
Let $G$ be a discrete group which admits an amenable action on a
compact space and $\gamma \in Aut(G)$ be an automorphism.  We define a
notion of entropy for $\gamma$ and denote the invariant by
$ha(\gamma)$.  This notion is dual to classical topological entropy in
the sense that if $G$ is abelian then $ha(\gamma) =
h_{Top}(\hat{\gamma})$ where $h_{Top}(\hat{\gamma})$ denotes the
topological entropy of the induced automorphism $\hat{\gamma}$ of the
(compact, abelian) dual group $\hat{G}$.

$ha(\cdot)$ enjoys a number of basic properties which are useful for
calculations.  For example, it decreases in invariant subgroups and
certain quotients.  These basic properties are used to compute the dual
entropy of an arbitrary automorphism of a crystallographic group.
\end{abstract}

\parskip2mm

\section{Introduction}

Let $\cal G$ denote the class of countable discrete groups which admit
an amenable action on a compact space (cf.\ [ADR]; see also Definition
2.1 below).  This notion has recently been very important in work
related to the Novikov conjecture (cf.\ [HR]).  Recently it has also
been shown that this class of groups coincides with the class of exact
groups (cf.\ [AD],[GK],[O]).  $\cal G$ is known to be very large; it
contains the classical amenable groups, free groups (or any other
hyperbolic group; [Ad], [Ge]) and discrete subgroups of connected Lie
groups (cf.\ [Z]).  Moreover, $\cal G$ is closed under a number of
standard group operations such as taking subgroups, extensions,
 free products  and quotients by classical amenable groups.

In this paper we study a conjugacy class invariant for automorphisms
of elements in $\cal G$.  The invariant, denoted $ha(\gamma)$ when
$\gamma \in Aut(G)$ and $G \in \cal G$, is called the {\em dual entropy} of
$\gamma$.  The terminology comes from the fact that if $G$ is abelian
with (compact) dual group $\hat{G}$ then for each $\gamma \in Aut(G)$,
$ha(\gamma)$ is equal to the classical topological entropy of the
induced automorphism $\hat{\gamma} \in Aut(\hat{G})$.   

Our motivation for the present work is twofold.  On the one hand,
discrete groups are one of the most basic objects in mathematics and
hence we feel it is worthwhile to have a numerical invariant for
dynamical systems arising from such groups.  In the framework of
operator algebras, noncommutative dynamical systems have been studied
for a long time.  However, the various notions of entropy in this
setting are still not well understood.  Hence the second motivation
for our work is to provide an invariant which is closely related to
some of the operator algebra entropies and will hopefully provide a
setting where hypotheses in general operator algebras can be more
easily tested.  For example, a result of Voiculescu concerning certain
automorphisms of noncommutative tori will follow immediately from a
general inequality relating dual entropy and noncommutative
topological entropy as defined in [Br] (cf.\ Propositions 3.3, 3.5).
Also, we will show that dual entropy decreases in quotients when the
kernel is a classical amenable group (cf.\ Proposition 5.8). (The
corresponding question for noncommutative topological entropy is still
open.)

These notes are organized as follows.  In section 2 we give the
definition of dual entropy and show that there is always a canonical
action which determines entropy calculations (Proposition 2.6).  In
section 3 we prove a general relationship between dual entropy and the
noncommutative topological entropy defined in [Br] (cf.\ Proposition
3.3).  In section 4 we justify the terminology "dual" entropy by
treating the case of abelian groups.  In section 5 we develop the
basic properties of this entropy function.  We will see that dual
entropy satisfies all of the properties which are currently known for
other noncommutative approximation entropies and has the added feature
of decreasing in certain quotients.  Finally in section 6 we give some
calculations of dual entropy.  In particular, we will give exact
calculations of the dual entropy of an arbitrary automorphism of a
crystallographic group (Corollary 6.7).

\section{Definition of Dual Entropy}

In this section, we define the dual entropy of an automorphism (or any
other map) of a discrete group which admits an amenable action.  Our
definition is based on the {\em approximation approach} to
noncommutative dynamical entropy which was first introduced in [Vo].
In [Vo] the main focus is operator algebras and the approximations
involved deal with (various) norm approximations of operators on a
Hilbert space.  In the group context we replace norm approximations
with {\em approximate invariance} and get a definition which behaves
similar to those in [Vo] (see also [Br]).  In a recent preprint,
M. Choda has independantly formulated a notion similar to our dual
entropy in the setting of amenable groups (cf.\ [Ch]).  It would be
interesting to know whether or not these notions coincide for amenable
groups.

Before we can state the appropriate definitions, we must introduce
some notation.  If $G$ is a discrete group and $A$ is a unital
commutative $C^*$-algebra we let $l^1 (G, A)$ denote the closure of the linear
space of finitely supported functions $T : G \to A$ with respect to
the norm $\| T \|_1 = \| \sum_g | T(g) | \|_{A}$.  

\begin{definition}
{\em Let $G$ be a discrete group and $\alpha : G \to Aut(A)$ be a
homomorphism, where $A$ is a unital {\em commutative} $C^*$-algebra.
  Then the action
$\alpha$ is called {\em amenable} if there exist functions $T_n \in
l^1 (G, A)$ such that $T_n$ is nonnegative (i.e. $T_n (g) \geq 0$ for
all $g \in G$), finitely supported, $\sum_g T_n (g) = 1_A$ and $\|
s.T_n - T_n \|_1 \to 0$ for all $s \in G$, where $s.T
(g) = \alpha_s (T(s^{-1}g))$ for all $T \in l^1 (G, A)$.}
\end{definition}

It is clear that a group is amenable in the classical sense if and
only if the trivial action $\tau : G \to Aut({\Bbb C})$ is amenable in
the sense described above.  See also [ADR] for a comprehensive
treatment of general amenable groupoids. 

\begin{definition}
{\em If $G$ is a discrete group, $\alpha : G \to Aut(A)$ is an
amenable action, $\omega \subset G$ is a finite set and $\delta > 0$
then put $$ ra(\alpha, \omega, \delta) = \inf\{ |supp(T)| : \| s.T -
T \|_1 < \delta \ for \ all \ s \in \omega\},$$ where the infimum is
taken over all nonnegative functions $T \in l^1(G,A)$ with $\sum_g
T(g) = 1_A$ and $|supp(T)|$ denotes the cardinality of the support of
$T$.  The integer $ra(\alpha, \omega, \delta)$ is called the {\em
amenable $\delta$-rank of $\omega$ with respect to $\alpha$}. One then
defines $$ ra(\omega, \delta) = \inf_{\alpha} ra(\alpha, \omega,
\delta),$$ where the infimum is taken over all amenable actions $\alpha
: G \to Aut(A)$ (and, of course, $A$ is allowed to vary).  Then
$ra(\omega, \delta)$ is called the {\em amenable $\delta$-rank of
$\omega$}.}
\end{definition}

We will soon see that the infimum defining $ra(\omega, \delta)$ is
always realized at a canonical action (Proposition 2.6).  However, the
present definition makes certain calculations easier.

Having a suitable ``$\delta$-rank'' function for finite sets we now
mimic the definitions in [Vo] to get our definition of dual entropy.

\begin{definition}
{\em If $G$ admits an amenable action and $\gamma \in Aut(G)$ we define the
following quantities: $$ha(\gamma, \omega, \delta) = \limsup_{n \to
\infty} \frac{\log\big( ra(\omega \cup \gamma(\omega) \cup \ldots \cup
\gamma^{n - 1}(\omega), \delta) \big)}{n},$$ $$ha(\gamma, \omega) =
\sup_{\delta > 0} \ ha(\gamma, \omega, \delta),$$ $$ha(\gamma) =
\sup_{\omega} \ ha(\gamma, \omega),$$ where the final supremum is
taken over all finite subsets of $G$.  We call $ha(\gamma)$ the {\em
dual entropy of $\gamma$}.}
\end{definition}

Though stated for automorphisms, it is clear that the definitions
above make perfectly good sense for endomorphisms (or any other self
map) of $G$.  Indeed, many of the results which are to follow have
similar formulations and proofs in the case of endomorphisms.

\begin{remark}
{\em The above definitions are also easily extended to cover arbitrary
(locally compact) groups which admit amenable actions.  Since Haar
measure on a discrete group is just the counting measure, we could
also write $\mu(supp(T))$ in the definition of $ra(\alpha, \omega,
\delta)$, where $\mu$ is Haar measure, and then it is clear how to
extend this definition to cover more general groups.  However, we will
stick to the discrete case since this more general notion would have
the philosophically unpleasant feature of always taking the value zero
on any compact group. }
\end{remark}

\begin{definition}
{\em If $\alpha : G \to Aut(A)$ is an amenable action, $\omega \subset
G$ is a finite set and $\delta > 0$ then we say a funtion $T \in l^1
(G, A)$ is {\em minimal for $(\omega, \delta)$} if $T$ is nonnegative,
$\sum_g T(g) = 1_A$, $\| s.T - T \|_1 < \delta$ for all $s \in
\omega$ and $|supp(T)| = ra(\alpha, \omega, \delta)$.}
\end{definition}

If $G$ is discrete we let $\alpha^G : G \to Aut(l^{\infty}(G))$ denote
the canonical action by left translation.  That is, $\alpha^G_g
(\zeta) = g.\zeta$ where $g.\zeta (x) = \zeta(g^{-1}x)$ for all
$\zeta \in l^{\infty}(G)$.  It is well known that $G$ admits an
amenable action if and only if the action $\alpha^G$ is already
amenable.  However, it is an important fact that all entropy
calculations are also determined by this canonical action.

\begin{proposition}
If $G$ admits an amenable action then for all finite subsets $\omega$
and $\delta > 0$, $ra(\omega, \delta) = ra(\alpha^G, \omega, \delta)$.
If $G$ is an amenable group then $ra(\omega, \delta) = ra(\tau,
\omega, \delta)$, where $\tau : G \to Aut({\Bbb C})$ is the trivial
action.
\end{proposition}

{\noindent\bf Proof.}  Let $\omega \subset G$ be a finite set and
$\delta > 0$.  Choose an amenable action $\alpha : G \to Aut(C(X))$
for which $ra(\omega, \delta) = ra(\alpha, \omega, \delta)$ and take $
T \in l^1 (G, C(X))$ minimal for $(\omega, \delta)$ (Definition 2.5).

 Fix some point $x_0 \in X$ and consider the
unital *-homomorphism $\Lambda : C(X) \to l^{\infty}(G)$ given by
$\xi\mapsto \Lambda(\xi)$ such that $\Lambda(\xi)(g) =
\alpha_g(\xi)(x_0)$
 for all $g \in G$ and
$\xi \in C(X)$.  The important observation is that $\Lambda$
intertwines the actions $\alpha$ and $\alpha^G$.

Recall that we have $T \in l^1 (G, C(X))$ minimal for $(\omega,
\delta)$ and hence we define $S \in l^1 (G, l^{\infty}(G))$ by $S(g) =
\Lambda(T(g))$.  Then it is clear that $S$ is nonnegative, $\sum_g
S(g) = 1_{l^{\infty}(G)}$ and $supp(S) \subset supp(T)$.  Moreover, for each
$h \in G$ we have 
\begin{multline*}
\begin{aligned}
\| h.S - S \|_1  
    &=  \| \sum_{g \in G} | \alpha^G_h(\Lambda(T(h^{-1}g)) ) 
           - \Lambda(T(g)) |  \ \|_{l^{\infty}(G)}  \\[2mm]
    &=  \| \Lambda\bigg( \sum_{g \in G} | \alpha_h( T(h^{-1}g)) 
           - T(g)) | \bigg)  \ \|_{l^{\infty}(G)}   \\[2mm]
    &\leq  \|  \sum_{g \in G} | \alpha_h( T(h^{-1}g)) 
           - T(g)) |  \ \|_{C(X)}   \\[2mm] 
    &=  \| h.T - T \|_1.
\end{aligned}
\end{multline*}
Thus we see that $ra(\alpha^G, \omega, \delta) \leq |supp(S)| \leq
|supp(T)| = ra(\omega, \delta)$.  However, the opposite inequality is
immediate from the definition and so we have proved the first part of
the proposition. 

Assume now that $G$ is an amenable group and let $m$ be a left
invariant mean on $l^{\infty}(G)$ (i.e. $m$ is a state on
$l^{\infty}(G)$ with the property that $m(\alpha^G_h(\zeta)) =
m(h.\zeta) = m(\zeta)$ for all $\zeta \in l^{\infty}(G)$ and $h \in
G$).  Let $\tau : G \to Aut({\Bbb C})$ be the trivial action and by
the first part of the proof we only have to show that $ra(\tau,
\omega, \delta) \leq ra(\alpha^G, \omega, \delta)$ for all finite
subsets $\omega \subset G$ and $\delta > 0$.  So choose $T \in l^1 (G,
l^{\infty}(G))$ minimal for $(\omega, \delta)$.  Then define $S \in
l^1(G) = l^1 (G, {\Bbb C})$ as $S(g) = m(T(g))$ and for all $h \in G$
we have
\begin{multline*}
\begin{aligned}
\| h.S - S \|_1  
    &=   \sum_{g \in G} | S(h^{-1}g) - S(g) | \\[2mm] 
    &=   \sum_{g \in G} | m(T(h^{-1}g)) - m(T(g)) |\\[2mm]
    &=   \sum_{g \in G} | m(\alpha^G_h(T(h^{-1}g))) - m(T(g)) |\\[2mm]
    &\leq m\bigg( \sum_{g \in G} | \alpha^G_h(T(h^{-1}g)) - T(g) | 
         \bigg) \\[2mm]
    &\leq \| h.T - T \|_1.
\end{aligned}
\end{multline*}
Since it is clear that $S$ is nonnegative, $\sum_g S(g) = 1$ and
$supp(S) \subset supp(T)$ (actually, we have equality by the
minimality of $T$ and the estimates above) we see that $ra(\tau,
\omega, \delta) \leq ra(\alpha^G, \omega, \delta) = ra(\omega, \delta)$.  
$\Box$

\section{Comparison with Noncommutative Topological Entropy}

In this section we show that dual entropy on a discrete group
dominates the noncommutative topological entropy (as defined in [Br])
on the reduced group $C^*$-algebra.  In fact, the same inequality
holds when the group has a $\Bbb C$ valued cocycle and hence we will
immediately deduce a result of Voiculescu for certain natural
automorphisms of noncommutative tori.  See [Br] for all the
definitions and notation related to noncommutative topological
entropy.

If $\alpha : G \to Aut(A)$ in an action on a unital commutative $C^*$-algebra and
$S, T : G \to A$ are finitely supported functions we define the
$\alpha$-convolution as follows: $$S \ast_{\alpha} T (h) = \sum_{g \in
G} S(g)\alpha_g (T(g^{-1}h)).$$ We also define the function $T^*$ by
$T^* (h) = \alpha_h(T(h^{-1})^*)$.  Let $l^2 (G, A)$ denote the
closure of the finitely supported functions from $G$ to $A$ with
respect to the norm $\| T \|_2 = \| \sum_g T(g)^* T(g) \|_{A}^{1/2}$.
$l^2 (G, A)$ is a Hilbert $A$-module with respect to the inner product
$<T,S> = \sum_g T(g)^* S(g)$ and we have the usual Cauchy-Schwartz
inequality $\| <T,S> \|_{A} \leq \| T \|_2 \| S \|_2$.

\begin{lemma}
If  $T : G \to A$ is finitely supported, nonnegative,
$\sum_g T(g) = 1_A$ and $S : G \to A$ is defined by $S(g) =
T(g)^{1/2}$, then for all $h \in G$ we have $\| 1_A -
S\ast_{\alpha}S^* (h) \|_{A}^{2} \leq \| h.T - T \|_{l^1(G, A)}$.
\end{lemma}

{\noindent\bf Proof.}  The proof is a calculation analogous to that in
[Pe, 7.3.8] only a bit more care must be taken since we are dealing
with modules and functions. 
\begin{multline*}
\begin{aligned}
\| 1_A - S\ast_{\alpha}S^* (h) \|_{A}^{2} 
    &=   \| <S, S - h.S> \|_{A}^{2}    \\[2mm]
    &\leq  \| S - h.S \|_2^2  \\[2mm]
    &=   \| \sum_{g \in G} | S(g) - \alpha_h(S(h^{-1}g)) |^2 \|_A \\[2mm]
    &\leq  \| \sum_{g \in G} | S(g)^2 - \alpha_h(S(h^{-1}g))^2 | \ \|_A \\[2mm]
    &= \| T - h.T \|_1 .   \ \ \ \ \ \Box
\end{aligned}
\end{multline*}

\begin{definition}
{\em A map $\theta : G \times G \to {\Bbb T} = \{ c \in {\Bbb C} :
|c| = 1 \}$ is called a ${\Bbb C}$-cocycle if $\theta(1,g) =
\theta(g,1) = 1$, for all $g \in G$ and $\theta(g,h)\theta(gh,k) =
\theta(h,k)\theta(g,hk)$ for all $g,h,k \in G$.  }
\end{definition}

If $\alpha : G \to Aut(A)$ is an action where $A \subset B(H)$ then we
define a *-monomorphism $\pi : A \to B(l^2 (G) \otimes H)$ by $\pi(a)
\xi_g \otimes h = \xi_g \otimes \alpha_{g^{-1}}(a)h$ for all $\xi_g
\in l^2 (G)$ and $h \in H$. ( $\{ \xi_g \}_{g \in G}$ denotes the
canonical orthonormal basis.)  We also define a unitary representation
$\lambda^{\theta}: G \to B(l^2 (G) \otimes H)$ by $\lambda^{\theta}_s
(\xi_g \otimes h) = \theta(s,g) \xi_{sg} \otimes h$.  Then the reduced
cocycle crossed product $A \rtimes_{r,\alpha} (G,\theta)$ is defined to
be the $C^*$-algebra generated by $\{ \pi(a)\lambda^{\theta}_s : a \in
A, s \in G\}$.  One has the relations
$\lambda^{\theta}_s\lambda^{\theta}_t = \theta(s,t)
\lambda^{\theta}_{st}$ and $\lambda^{\theta}_s \pi(a)
\lambda^{\theta\ast}_s = \pi(\alpha_s(a))$.  When $A = {\Bbb C}$ we
use the notation $C_r^* (G,\theta)$.  Note that if $\gamma \in Aut(G)$
and satisfies the relation $\theta(g,h) = \theta(\gamma(g),
\gamma(h))$ for all $g,h \in G$ then there is a natural induced
automorphism $\hat{\gamma} \in Aut(C_r^* (G,\theta))$ defined by
$\hat{\gamma}( \lambda^{\theta}_g) = \lambda^{\theta}_{\gamma(g)}$.

\begin{proposition}
If $\theta : G\times G \to {\Bbb C}$ is a $\Bbb C$-cocycle and $\gamma
\in Aut(G)$ is such that $\theta(g,h) = \theta(\gamma(g), \gamma(h))$
for all $g,h \in G$ then $ht(\hat{\gamma}) \leq ha(\gamma)$.
\end{proposition}

{\noindent\bf Proof.} The proof is an adaptation of the techniques in
[SS], [Br] and [BC]: thus we will be rather sketchy.  Let $\alpha^G
: G \to Aut(l^{\infty}(G))$ be the canonical action and construct
$l^{\infty}(G) \rtimes_{r,\alpha^G} (G,\theta)$ as above.

By the Kolmogorov-Sinai type result ([Br, Prop.\ 2.6]) and
Proposition 2.6 it suffices to show that if $\omega \subset G$ is a
finite set then $rcp(\iota, \lambda^{\theta}(\omega), \delta^{1/2})
\leq ra(\alpha^G, \omega, \delta)$ for all $\delta > 0$, where $\iota
: l^{\infty}(G) \rtimes_{r,\alpha^G} (G,\theta) \hookrightarrow B(l^2
(G) \otimes H)$ is the natural inclusion. So, choose $T \in l^1(G,
l^{\infty}(G))$ which is minimal for $(\omega, \delta)$ and let $F =
supp(T) \subset G$.  Let $P_F \in B(l^2(G))$ be the orthogonal
projection onto the span of $\{ \xi_h : h \in F \}$.  Let $X = \sum_{p
\in F} e_{p,p} \otimes \alpha^{G}_{p^{-1}}(S(p)) \in M_{|F|}\otimes
l^{\infty}(G)$, where $S(p) = T(p)^{1/2}$ and $e_{p,q}$ denote the
canonical matrix units of $M_{|F|} = P_F B(l^2(G))P_F$.  Then one
computes $$X(P_F \otimes I)\lambda^{\theta}_h(P_F \otimes I)X =
\sum_{p \in F \cap hF} \theta(h,h^{-1}p)e_{p, h^{-1}p}\otimes
\alpha^{G}_{p^{-1}}(S(p))\alpha^{G}_{p^{-1}h}(S(h^{-1}p)).$$ Now
define the completely positive map $L : M_{|F|}\otimes l^{\infty}(G)
\to l^{\infty}(G) \rtimes_{r,\alpha^G} (G,\theta)$ by $L (e_{x,y}
\otimes \zeta) = \overline{\theta(xy^{-1}, y)} \pi(\alpha^G_x
(\zeta))\lambda^{\theta}_{xy^{-1}}$.  Then a straightforward (but
somewhat unpleasant) calculation shows that $$L(X(P_F \otimes
I)\lambda^{\theta}_h(P_F \otimes I)X) =
\pi(S\ast_{\alpha^G}S^*(h))\lambda_h^{\theta},$$ for all $h \in
G$. Hence $\| \lambda_h^{\theta} - L(X(P_F \otimes
I)\lambda^{\theta}_h(P_F \otimes I)X) \| \leq \delta^{1/2}$ for all $h
\in \omega$, by Lemma 3.1. But since the map $\lambda_h^{\theta}
\mapsto L(X(P_F \otimes I)\lambda^{\theta}_h(P_F \otimes I)X)$ factors
through the matrix algebra $M_{|F|}$ this implies the desired
inequality.  $\Box$

\begin{remark} 
{\em It follows from [Dy, Prop.\ 9] and the proposition above that
the dual entropy also dominates CNT-entropy (cf.\ [CNT]) with respect
to the canonical trace on $C^*_r (G)$ (and in the $W^*$-algebra
$L(G)$).  Also, the proof above together with the results of [BC,
Section 2]  can be used to show that dual entropy dominates the
$W^*$-entropy defined in [Vo, Section 3] in the case that $G$ is an
amenable group. }
\end{remark}

Let $\theta : {\Bbb Z}^n \times {\Bbb Z}^n \to {\Bbb T}$ be a $\Bbb
C$-cocycle.  Then $C_r^* ({\Bbb Z}^n ,\theta)$ is the usual
noncommutative torus with generating unitaries $\lambda^{\theta}_1,
\ldots, \lambda^{\theta}_n$ ($\lambda^{\theta}_i =
\lambda^{\theta}_{{e_i}}$ for the canonical generators $\{ e_i \}_{1
\leq i \leq n}$ of ${\Bbb Z}^n$) subject to the relations $
\lambda^{\theta}_i \lambda^{\theta}_j = \theta(e_i,e_j)
\lambda^{\theta}_{i + j} = \theta(e_i,e_j)
\overline{\theta(e_j,e_i)}\lambda^{\theta}_j\lambda^{\theta}_i$. 

\begin{corollary}[Vo, Prop.\ 5.3]
If $\gamma \in GL(n, {\Bbb Z}^n) = Aut({\Bbb Z}^n)$ is a matrix such
that $\theta(\gamma(g), \gamma(h)) = \theta(g,h)$ for all $g,h \in
{\Bbb Z}^n$ then $ht(\hat{\gamma}) \leq ha(\gamma) = \log(\mu_1 \mu_2
\cdots \mu_n)$, where $\hat{\gamma}$ is the induced automorphism of
the noncommutative torus $C_r^* ({\Bbb Z}^n ,\theta)$ and $\mu_j =
\max(1, |t_j|)$, $1 \leq j \leq n$, for the eigenvalues
$t_j$ of $\gamma$.
\end{corollary}

{\noindent\bf Proof.}  In the next section we will show that
$ha(\gamma) = \log(\mu_1 \mu_2 \cdots \mu_n)$, since the dual group of
${\Bbb Z}^n$ is the (commutative) $n$-torus and it is well known how
to compute the topological entropy of automorphisms of ${\Bbb T}^n$.
Hence the corollary follows from Proposition 3.3.  $\Box$

\section{The Abelian Case}

Our next goal is to justify the terminology ``dual'' entropy.  We will
need a theorem of J. Peters (see [Pet, Thm.\ 6]).

\begin{theorem}
Let $G$ be a discrete abelian group with dual group $\hat{G}$. If
$\gamma \in Aut(G)$ with induced automorphism $\hat{\gamma} \in
Aut(\hat{G})$ then the (classical) topological entropy,
$h_{Top}(\hat{\gamma})$, of $\hat{\gamma}$ is equal to $$\sup_{E}
\bigg( \limsup_{n \to \infty}
\frac{\log |E + \gamma(E) + \ldots + \gamma^{n - 1}(E)|}{n}\bigg),$$
where the supremum is taken over all finite subsets $E \subset G$.
\end{theorem}

\begin{theorem}
If $G$ is discrete and abelian, $\gamma \in Aut(G)$ and $\hat{\gamma}
\in Aut(\hat{G})$ is the induced automorphism then $ha(\gamma) =
h_{Top}(\hat{\gamma})$.
\end{theorem}

{\noindent\bf Proof.}  By Proposition 3.3 above, [Br, Prop.\ 1.4] and
[Vo, Prop.\ 4.8] we have the inequality $ha(\gamma) \geq
h_{Top}(\hat{\gamma})$.  To prove the other inequality we apply
Proposition 2.6 (recall that abelian groups are amenable) and Peters'
theorem above.  One basic fact we will need is that if $f_1, f_2 \in
l^1 (G)$ are nonnegative functions then $ \| f_1 \ast f_2 \|_1 =
\| f_1 \|_1 \| f_2 \|_1$ and $supp( f_1 \ast f_2 ) = supp(f_1) +
supp(f_2)$.

So let $\omega \subset G$ be a finite set and $\delta > 0$.  Choose
some nonnegative, norm one function of finite support $f \in l^1 (G)$
such that $\| s.f - f\| < \delta$ for all $s \in \omega$.  Now define
$F_n = f \ast f\circ \gamma^{-1} \ast \cdots \ast f\circ \gamma^{-n +
1}$.  Then one verifies (using that convolution is commutative since
$G$ is abelian) that $\| s.F_n - F_n \| < \delta$ for all $s \in
\omega \cup \cdots \cup \gamma^{n - 1}(\omega)$.  Hence $ra(\omega
\cup \cdots \cup \gamma^{n - 1}(\omega), \delta) \leq |supp(F_n)|$.
But $supp(F_n) = supp(f) + \cdots + \gamma^{n - 1}(supp(f))$ and hence
one deduces that $ha(\gamma, \omega, \delta)$ is less than or equal to
Peters' formula and hence is bounded above by $h_{Top}(\hat{\gamma})$.
$\Box$

\section{Basic Properties}

In this section we develop the basic properties of the fuction $ha(
\cdot )$.  Some of the proofs run parallel to the corresponding
properties for the approximation entropies defined in [Vo] and hence
we will often refer to that paper for details.

We begin by observing that $ra( \cdot, \cdot)$ is an isomorphism
invariant.

\begin{lemma} 
Let $\Phi : G \to H$ be a group isomorphism.  For all
finite subsets $\omega \subset G$ and $\delta > 0$ we have $ra(\omega,
\delta) = ra(\Phi(\omega), \delta)$.  
\end{lemma}

{\noindent\bf Proof.}  Let $\alpha : G \to Aut(A)$ be an amenable
action and take $T \in l^1 (G, A)$ minimal for $(\omega, \delta)$
(cf.\ Definition 2.5).  Define an action $\beta : H \to Aut(A)$ by
$\beta(h) = \alpha_{\Phi^{-1}(h)}$ and $S \in l^1 (H, A)$ by $S(h) =
T(\Phi^{-1}(h))$.  Then one checks that $\| h.S - S \|_1 = \| g.T - T
\|_1$ for all $h = \Phi(g) \in H$. This evidently implies $ra(\omega,
\delta) \geq ra(\Phi(\omega), \delta)$. The opposite inequality is
similar. $\Box$

\begin{proposition} 
For all $\gamma, \sigma \in Aut(G)$, $ha(\gamma) = ha(\sigma \circ
\gamma \circ \sigma^{-1})$.
\end{proposition}

{\noindent\bf Proof.}  This is an easy consequence of the previous
lemma. $\Box$

\begin{proposition} 
For all $\gamma \in Aut(G)$ and $k \in {\Bbb Z}$ we have $ha(\gamma^k)
= |k| ha(\gamma)$.
\end{proposition}

{\noindent\bf Proof.}  With Lemma 5.1 in hand, the proof of this is
quite similar to the proof of [Vo, Prop.\ 1.3] as it depends only on
the algebraic part of the definition of $ha(\cdot)$ and not on the
particular $\delta$-rank function being used.  $\Box$

\begin{proposition} 
If $\omega_1 \subset \omega_2 \subset \ldots$ are finite sets with the
property that $$G = \bigcup_{i \in {\Bbb N}, n \in {\Bbb Z}} \gamma^n
(\omega_i)$$ then $ha(\gamma) = \sup_{i \in {\Bbb N}} ha(\gamma,
\omega_i)$.
\end{proposition}

{\noindent\bf Proof.}  This is similar to the proof of [Vo, Prop.\
3.4].  However, the reader will likely find it easier to give a proof
independantly as one has no topological considerations in this
setting. $\Box$

\begin{proposition}[Monotonicity] 
If $\gamma \in Aut(G)$ and $H \subset G$ is a subgroup such that
$\gamma(H) = H$ then $ha(\gamma|_{H}) \leq ha(\gamma)$.
\end{proposition}

{\noindent\bf Proof.}  Let $A$ be abelian, $\alpha : G \to Aut(A)$ be
an amenable action and $\beta : H \to Aut(A)$ denote the restriction
of $\alpha$ to $H$.  It suffices to show that for all finite subsets
$\omega \subset H$ and $\delta > 0$, $ra(\beta, \omega, \delta) \leq
ra(\alpha, \omega, \delta)$.

Let $G/H$ be the space of right cosets and choose representatives $\{
p_i \}_{i \in G/H} \subset G$.  Then each element of $G$ has a unique
representation as $hp_j$ for some $h \in H$ and $j \in G/H$.  Now
choose $T \in l^1(G, A)$ which is minimal for $(\omega, \delta)$ and
define $S \in l^1(H, A)$ by $$S(h) = \sum_{i \in G/H} T(hp_i).$$ Since
$T$ is finitely supported, the summation above is finite.  It is clear
that $S$ is nonnegative, $\sum_{h \in H} S(h) = \sum_h \sum_i T(hp_i)
= 1_A$ and $|supp(S)| \leq |supp(T)|$.  Moreover, for all $h \in H$,
\begin{multline*} 
\begin{aligned} 
\| h.S - S \|_1
    &=   \| \sum_{k \in H} | \beta_{h}(S(h^{-1}k))  - 
         S(k) | \ \|  \\[2mm] 
    &=   \| \sum_{k \in H} | \sum_{i \in G/H} 
          \alpha_h (T(h^{-1}kp_i)) - T(kp_i) | \ \|  \\[2mm]
    &\leq \| \sum_{g \in G} | \alpha_h (T(h^{-1}g)) - T(g) | \ \|\\[2mm]
    &=   \| h.T - T \|_1.
\end{aligned}
\end{multline*}
This implies $ra(\beta, \omega, \delta) \leq ra(\alpha, 
\omega, \delta)$.   $\Box$

We now turn to the question of how $ha(\cdot)$ behaves in quotients.
We rather doubt that the dual entropy always decreases but it does
behave well when the kernel is amenable. However, this will require
some preliminary results.

So assume that $1 \to K \to G \stackrel{\pi}{\to} H \to 1$ is an exact
sequence where $K$ is an amenable group with left invariant mean $m
\in S(l^{\infty} (K))$.  Let $\phi : H \to G$ be a unital (set
theoretic) splitting (i.e. $\phi(1_H) = 1_G$ and $\pi(\phi(h)) = h$
for all $h \in H$) and $\theta : H \times H \to K$ the corresponding
cocycle map (i.e. $\phi(h_1)\phi(h_2) = \phi(h_1 h_2)\theta(h_1, h_2)$
for all $h_1, h_2 \in H$).  For each $\zeta \in l^{\infty} (G)$ and $h
\in H$, let $\zeta(\phi(h) \cdot) \in l^{\infty} (K)$ denote the
function $k \mapsto \zeta(\phi(h)k)$. Define a linear map $\Lambda :
l^{\infty} (G) \to l^{\infty} (G)$ by $$\Lambda(\zeta)(\phi(h)k) =
m(\zeta(\phi(h) \cdot))$$ for all $g = \phi(h)k \in G$.

\begin{lemma}
$\Lambda : l^{\infty} (G) \to l^{\infty} (G)$ is independant of the
splitting $\phi$.  Also, for all $s \in G$ and $\zeta \in l^{\infty}
(G)$, $\Lambda(\alpha^G_s(\zeta)) = \alpha^G_s(\Lambda(\zeta))$.
\end{lemma}

{\noindent\bf Proof.}  Let $\psi : H \to G$ be any other unital
splitting.  Then for each $h \in H$ there exists $k_h \in K$ such that
$\psi(h) = \phi(h)k_h$.  So for $\zeta \in l^{\infty} (G)$ we have
$\zeta(\psi(h) \cdot)(k) = \zeta(\phi(h) k_h k) =
\alpha^K_{k_h^{-1}}(\zeta(\phi(h) \cdot))(k)$.  Hence $m(\zeta(\psi(h)
\cdot)) = m(\zeta(\phi(h) \cdot))$ and we see that $\Lambda$ does not
depend on $\phi$. 

Now let $s = p\phi(q)^{-1} \in G$ be arbitrary.  For each $\zeta \in
l^{\infty} (G)$, $x \in H$ and $k \in K$ we have
\begin{multline*} 
\begin{aligned} 
\alpha^G_s(\zeta)(\phi(x) \cdot)(k)
    &=   \alpha^G_s(\zeta)(\phi(x)k) \\[2mm] 
    &=   \zeta(\phi(q)p^{-1}\phi(x)k) \\[2mm]
    &=   \zeta(\phi(qx)\theta(q,x)\phi(x)^{-1}p^{-1}\phi(x)k) \\[2mm]
    &=   \alpha^K_{t^{-1}}(\zeta(\phi(qx)\cdot))(k), 
\end{aligned}
\end{multline*}
where $t = \theta(q,x)\phi(x)^{-1}p^{-1}\phi(x)$ (note that
$\phi(x)^{-1}p^{-1}\phi(x) \in K$ by normality).  Hence
$m(\alpha^G_s(\zeta)(\phi(x) \cdot)) = m(\zeta(\phi(qx)\cdot))$.  This
implies that $\Lambda(\alpha^G_s(\zeta))(\phi(x)y) = m( \zeta(\phi(qx)
\cdot) )$.

A similar sort of argument shows $\alpha^G_s(\Lambda(\zeta))(\phi(x)y)
= m( \zeta(\phi(qx) \cdot) )$ and hence $\Lambda(\alpha^G_s(\zeta)) =
\alpha^G_s(\Lambda(\zeta))$ as claimed.  $\Box$

\begin{lemma}
If $\alpha^G : G \to Aut(l^{\infty} (G))$ is amenable and $T \in
l^1(G, l^{\infty} (G))$ is minimal for $(\omega, \delta)$ then
$\Lambda(T) \in l^1(G, l^{\infty} (G))$ is also minimal for $(\omega,
\delta)$ where $\Lambda(T)(g) = \Lambda(T(g))$.
\end{lemma}

{\noindent\bf Proof.} Evidently $\Lambda$ is a positive linear map
(i.e. if $\zeta \geq 0$ then $\Lambda(\zeta) \geq 0$).  Hence if
$\zeta \in l^{\infty} (G)$ takes values in $\Bbb R$ then
$|\Lambda(\zeta)| \leq \Lambda(|\zeta|)$ since $\zeta \leq |\zeta|$.
With this observation and the above fact that
$\Lambda(\alpha^G_s(\zeta)) = \alpha^G_s(\Lambda(\zeta))$ for all
$\zeta \in l^{\infty} (G)$, a routine calculation shows that $\|
g.\Lambda(T) - \Lambda(T) \|_1 \leq \| g.T - T \|_1$ for all $g
\in G$.  This implies the lemma since it is clear that $\Lambda(T)$ is
nonnegative, $supp(\Lambda(T)) \subset supp(T)$ (actually $=$ by the
minimality of $T$) and $\sum_g \Lambda(T)(g) = 1_A$.  $\Box$

\begin{proposition} 
Assume $1 \to K \to G \stackrel{\pi}{\to} H \to 1$ is an exact
sequence, where $K$ is an amenable group, and assume $\gamma \in
Aut(G)$ leaves $K$ invariant (i.e. $\gamma(K) = K$).  If $\dot{\gamma}
\in Aut(G/K)$ denotes the induced automorphism then $ha(\gamma) \geq
ha(\dot{\gamma})$.
\end{proposition}

{\noindent\bf Proof.}  It suffices to show that for each finite subset
$\omega \subset G$ and $\delta > 0$, $ra(\omega, \delta) \geq
ra(\pi(\omega), \delta)$. 

Let $\phi : H \to G$ be a unital splitting with cocycle map $\theta :
H \times H \to K$ and construct $\Lambda : l^{\infty} (G) \to
l^{\infty} (G)$ as above.  Let $\pi_* : l^{\infty} (H) \hookrightarrow
l^{\infty} (G)$ denote the unital *-monomorphism induced by $\pi : G
\to H$.  One readily verifies that $\pi_* (\alpha^H_{\pi(x)} (\zeta))
= \alpha^G_x (\pi_*(\zeta))$ for all $x \in G, \zeta \in l^{\infty}
(H)$. Note also that $\pi_*(l^{\infty} (H)) = \Lambda(l^{\infty} (G))$
as subsets of $l^{\infty} (G)$.

Now choose $T \in l^1(G, l^{\infty} (G))$ minimal for $(\omega,
\delta)$.  By Lemma 5.7 above we may assume that $T(g) \in
\Lambda(l^{\infty} (G)) = \pi_*(l^{\infty} (H))$ for all $g \in G$.
Now define $S \in l^1(H, l^{\infty} (H))$ by $$S(\pi(x)) = \sum_{k \in
K} \pi_{*}^{-1}(T(xk))$$ for all $h = \pi(x) \in H$.  Note that
$S$ is nonnegative, $\sum_h S(h) = 1_{l^{\infty} (H)}$ and $|supp(S)|
\leq |supp(T)|$.  Finally we compute
\begin{multline*} 
\begin{aligned} 
\| \pi(x).S - S \|_1 
    &=  \| \sum_{h \in H} | \alpha^H_{\pi(x)} (S(\pi(x)^{-1}h)) 
        - S(h)| \ \|_{l^{\infty} (H)} \\[2mm]
    &=  \| \! \sum_{h \in H} \! \big| \! \sum_{k \in K} \!
        \pi_{*}^{-1}\bigg(\alpha^G_x(T(x^{-1}\phi(h)k)) - T(\phi(h)k) 
        \bigg) \big| \ \|_{l^{\infty} (H)} \\[2mm]
    &=  \| \sum_{h \in H} \big|  \sum_{k \in K} \alpha^G_x 
        (T(x^{-1}\phi(h)k)) - T(x^{-1}\phi(h)k) \big| \ 
        \|_{l^{\infty} (G)} \\[2mm]
    &\leq \| \sum_{g \in G} | \alpha^G_x(T(x^{-1}g)) - T(x^{-1}g) | \  
        \|_{l^{\infty} (G)}     \\[2mm]
    &=  \| x.T - T \|_1 .  \\[2mm]
\end{aligned}
\end{multline*}
This together with Proposition 2.6 implies that $ra(\omega, \delta) \geq
ra(\pi(\omega), \delta)$.  $\Box$

The next lemma is known to the experts so we only sketch the proof.

\begin{lemma}
If $H \subset G$ is a subgroup which admits an
amenable action, $\omega \subset H$ is a finite set and $\delta > 0$
then there exists $S \in l^1(G, l^{\infty} (G))$ such that $S$ is
nonnegative, has finite support, $\sum_g S(g) = 1_{l^{\infty} (G)}$ and
$\|h.S - S \| < \delta$ (w.r.t. the canonical action $\alpha^G$) for
all $h \in \omega$.
\end{lemma}

{\noindent\bf Proof.}  As in Proposition 5.5 we write $G = \bigcup_{k
\in H, i \in G/H} kp_i$ where $G/H$ denotes the space of right cosets.
Define a unital *-monomorphism $\Lambda : l^{\infty} (H) \to
l^{\infty} (G)$ by $\Lambda(\zeta)(kp_i) = \zeta(k)$.  A routine
calculation shows that for all $h \in H$, $\alpha^G_h (\Lambda(\zeta))
= \Lambda(\alpha^H_h(\zeta))$.  Taking $T \in l^1(H, l^{\infty} (H))$
minimal for $(\omega, \delta)$ and defining $S(g) = \Lambda (T(g))$ if
$g \in H$ and $S(g) = 0$ otherwise, one checks that $S \in l^1(G,
l^{\infty} (G))$ has the desired properties.  $\Box$

Let $G_1 \stackrel{\phi_1}{\to} G_2 \stackrel{\phi_2}{\to} G_3
\stackrel{\phi_3}{\to} \cdots$ be an inductive system of groups and
assume that there exist automorphisms $\gamma_i \in Aut(G_i)$ such
that $\gamma_{i + 1} \circ \phi_i = \phi_i \circ \gamma_i$ for all $i
\in {\Bbb N}$.  If $G = \lim\limits_{\longrightarrow} G_i$ is the
inductive limit then there is a natural inductive limit automorphism
$\gamma \in Aut(G)$ with the property that $\gamma \circ \Phi_i =
\Phi_i \circ \gamma_i$ for all $i \in {\Bbb N}$, where $\Phi_i : G_i
\to G$ are the natural homomorphisms.  Let $K_i = ker(\Phi_i) \subset
G_i$ and $\dot{\gamma_i}$ denote the induced automorphism of $G_i /
K_i$ for all $i \in {\Bbb N}$.

\begin{proposition} 
Let $G$ be an inductive limit of groups admitting amenable actions and
$\gamma \in Aut(G)$ be an inductive limit automorphism as above.  If
all of the subgroups $K_i$ are amenable (e.g. if the $\phi_i$ are
injective) then $$ha(\gamma) = \lim_{i \to \infty} ha(\dot{\gamma_i})
\leq \liminf_{i \to \infty} ha(\gamma_i).$$
\end{proposition}

{\noindent\bf Proof.}  Since Proposition 5.8 implies that the
groups $G_i / K_i$ admit amenable actions, it follows from Lemma 5.9
that $G$ also admits an amenable action.  Hence the proposition
follows from Propositions 5.4, 5.5 and 5.8.  $\Box$

We next investigate how $ha(\cdot)$ behaves in extensions.  Any
general results in this direction appear to be very hard as it is not
clear how to construct good actions on an extension when one is given
actions on the kernel and quotient.  The following lemma illustrates
the difficulties. Our proof is a modification of an argument shown to us
by Jean Renault and we thank him for sharing his notes with us.

As before, let $1 \to K \to G \stackrel{\pi}{\to} H \to 1$ be exact
and let $\phi : H \to G$ be a unital splitting with cocyle map $\theta
: H \times H \to K$.

\begin{lemma}
Let $G$ be as above and assume both $K$ and $H$ admit amenable
actions.  Let $\omega = \{ (\phi(h_1)k_1)^{-1}, \ldots,
(\phi(h_m)k_m)^{-1} \} \subset G$ be a finite set and choose $T \in
l^1(H, l^{\infty} (H))$ minimal for $(\pi(\omega), \delta)$.  If $F =
supp(T) \subset H$ then $$ra(\omega, 2\delta) \leq ra(\pi(\omega),
\delta)ra( \bigcup_{i = 1}^{m} \bigcup_{x \in h_i^{-1}F}
\phi(x)^{-1} k_i^{-1} \phi(x) \theta(h_i,x)^{-1}, \delta).$$
\end{lemma}

{\noindent\bf Proof.}  Choose $S \in l^1(K, l^{\infty} (K))$ minimal
for $$(\bigcup_{i = 1}^{m} \bigcup_{x \in h_i^{-1}F} \phi(x)^{-1}
k_i^{-1} \phi(x) \theta(h_i,x)^{-1}, \delta).$$ Now define $P : G \to
l^{\infty} (G)$ by $$P(\phi(x)y)(\phi(q)p) =
T(x)(q)S(y)(\phi(x)^{-1}\phi(q)p\phi(q^{-1}x)).$$ Note that
$\phi(x)^{-1}\phi(q)p\phi(q^{-1}x) =
(\phi(x)^{-1}(\phi(q)p\phi(q)^{-1})\phi(x))\theta(q,q^{-1} x) \in K$, by
normality of $K$, and hence $P$ is well defined.  It is clear that $P$
is nonnegative and finitely supported with $|supp(P)| \leq
|supp(T)||supp(S)|$.  Moreover, for all $\phi(q)p \in G$, 
\begin{multline*} 
\begin{aligned} 
\sum_{\phi(x)y \in G} P(\phi(x)y)(\phi(q)p)
    &=  \sum_{\phi(x)y \in G}  
        T(x)(q)S(y)(\phi(x)^{-1}\phi(q)p\phi(q^{-1}x))  \\[2mm]
    &=  \sum_{x \in H} T(x)(q) \big( \sum_{y \in K} 
         S(y)(\phi(x)^{-1}\phi(q)p\phi(q^{-1}x)) \big) \\[2mm]
    &=  \sum_{x \in H} T(x)(q) = 1.
\end{aligned}
\end{multline*}
Hence $\sum_g P(g) = 1_{l^{\infty} (G)}$.

It remains to show the inequality $\| s.P - P \|_1 < 2\delta$
for all $s = (\phi(h)k)^{-1} \in \omega$.  For notational reasons it
will be convenient to define $A(x) = \theta(h,x) \phi(x)^{-1} k
\phi(x)$ for all $x \in H$.
\begin{multline*} 
\begin{aligned} 
\| s.P - P \|_1
    &= \| \sum_{\phi(x)y \in G} | \alpha^G_{s}(P(\phi(h)k \phi(x)y)) 
       - P(\phi(x)y) | \ \|_{l^{\infty} (G)}  \\[2mm]
    &= \sup_{\phi(q)p} \bigg( \sum_{\phi(x)y \in G} | 
       P(\phi(hx)A(x)y) ( \phi(hq)A(q)p) \\[2mm]
    &\mathrel{\phantom\sup} \ \ - P(\phi(x)y)(\phi(q)p) | \bigg)  \\[2mm]
    &= \sup_{\phi(q)p} \!\! \bigg( \!\! \sum_{\phi(x)y \in G} \!\!\!| 
       T(hx)(hq) S(A(x)y)(A(x)\phi(x)^{-1}\phi(q)p\phi(q^{-1}x))   \\[2mm]
    &\mathrel{\phantom\sup} \ \ - T(x)(q)S(y)(\phi(x)^{-1}\phi(q)p
       \phi(q^{-1}x))  | \bigg)  \\[2mm]
    &\leq  \sup_{\phi(q)p} \!\!\bigg( \sum_{x \in h^{-1}F} \!\!T(hx)(hq) \big( 
           \sum_{y \in K} | S(A(x)y)(A(x)l) - S(y)(l) | \big) \\[2mm]
    &\mathrel{\phantom\sup} \ \ + \sum_{x \in H} |T(hx)(hq) - T(x)(q) | 
       \big( \sum_{y \in K} S(y)(l) \big) \bigg),   
\end{aligned}
\end{multline*}
where $l = \phi(x)^{-1}\phi(q)p\phi(q^{-1}x)$ in the last lines of the
inequality.  The inequality above implies the lemma.  $\Box$

\begin{corollary} 
Assume $1 \to K \to G \stackrel{\pi}{\to} H \to 1$
is exact with $K$ a finite group and $H$ admitting an amenable
action. If $\gamma \in Aut(G)$ is an automorphism such that $\gamma(K)
= K$ then $ha(\gamma) = ha(\dot{\gamma})$ where $\dot{\gamma} \in
Aut(H)$ is the induced automorphism.  
\end{corollary}

{\noindent\bf Proof.}  This is an easy consequence of Proposition 5.8
 and the previous lemma.  $\Box$

\begin{question}
In the previous corollary if we had required that $H$ be finite
(rather than $K$) would we then get $ha(\gamma) = ha(\gamma|_K)$?  We
are only able to solve this problem in the case that $K$ is a finitely
generated abelian group (cf.\ Theorem 6.4).
\end{question}

We now deduce the group analogues of [Br, Thm.\ 3.5] and [BC, Thm.\
3.3,3.4].  In particular, this gives an affirmative answer to the
analogue of [St, Problem 4.2].

\begin{proposition}
Let $K, \ H$ be groups admitting amenable actions and $\rho : H \to
Aut(K)$ be an action of $H$ on $K$.  If $h$ is in the center of $H$
then $ha(\rho_h) = ha({\rm Ad}h)$, where ${\rm Ad}h$ is the canonical
inner automorphism of $K \rtimes_{\rho} H$ implemented by $h$.
\end{proposition}

{\noindent\bf Proof.} Write $K \rtimes_{\rho} H = K \times H$ as sets
with multiplication $(k_1, h_1)(k_2, h_2) = (k_1
\rho_{h_1}(k_2), h_1h_2)$.  Let $\omega = \{ (\rho_{h_1}(k_1),
h_1)^{-1}, \ldots, (\rho_{h_m}(k_m), h_m)^{-1} \}$ be an arbitrary
finite set.  We also let $\pi : K \rtimes_{\rho} H \to H$ be the
canonical quotient map.  Note that for all $n \in {\Bbb N}$,
$\pi(\omega) = \pi( \omega \cup \ldots \cup {\rm Ad}h^n (\omega))$
since $h$ is central.

Now choose $T \in l^1(H, l^{\infty} (H))$ minimal for $(\pi(\omega),
\delta)$, let $F = supp(T)$ and consider the set $$ \tilde{\omega} =
\bigcup_{i = 1}^{m} \bigcup_{x \in h_i^{-1}F} (1_K,x)^{-1} (k_i,1_H)^{-1}
(1_K,x).$$ Using again the fact that $h$ is central, it follows from
Lemma 5.11 that $ra(\omega \cup \ldots \cup {\rm Ad}h^n (\omega),
2\delta)$ is bounded above by $$ra(\pi(\omega),
\delta)ra(\tilde{\omega} \cup \ldots \cup \rho_{h}^n (\tilde{\omega}),
\delta),$$ for all $n \in {\Bbb N}$.

Hence we deduce that $ha({\rm Ad}h, \omega, 2\delta) \leq
ha(\rho_h, \tilde{\omega}, \delta)$, which implies that $ha({\rm Ad}h)
\leq ha(\rho_h)$.  The opposite inequality follows from monotonicity
(Proposition 5.5).  $\Box$

\begin{corollary}
If $\gamma \in Aut(G)$ and $u \in G \rtimes_{\gamma} {\Bbb Z}$ is the
element which implements $\gamma$ in the semidirect product then
$ha(\gamma) = ha({\rm Ad}u)$.
\end{corollary}

\begin{proposition} 
For $\gamma_i \in Aut(G_i)$, $i = 1,2$, let $\gamma_1 \times \gamma_2$
denote the product automorphism of $G_1 \times G_2$.  Then $$\max \{
ha(\gamma_1), ha(\gamma_2) \} \leq ha(\gamma_1 \times \gamma_2) \leq
ha(\gamma_1) + ha(\gamma_2).$$
\end{proposition}

{\noindent\bf Proof.}  The lower bound follows from monotonicity. The
upper bound is a straightforward application of Lemma 5.11.  $\Box$

\begin{question}
{\em In the situation above do we have $ha(\gamma_1 \times \gamma_2) =
ha(\gamma_1) + ha(\gamma_2)$?  If the groups $G_1, G_2$ are abelian
then this equality holds since this is a well known property of
classical topological entropy (cf.\ Theorem 4.2). }
\end{question}

\begin{example} 
Let $G = \coprod_{\Bbb Z} G_0$ be the restricted
direct product of a finite group $G_0$ and let $\gamma \in Aut(G)$ be
the shift automorphism (induced by the mapping $i \mapsto i + 1$ of
$\Bbb Z$).  It follows from Proposition 3.3 that $\log(rank(C^*_r
(G_0))) \leq ha(\gamma)$, where $rank(M_{n_1}({\Bbb C}) \oplus \cdots
\oplus M_{n_k}({\Bbb C})) = n_1 + \cdots + n_k$, while the upper bound
$ha(\gamma) \leq \log(|G_0|)$ is readily seen.  It would be
interesting to know if we actually have equality at either the upper
or lower bound.  The difficulty appears to be the same as that
encountered in Question 5.17.  
\end{example}

It is known that if $G_1, \ G_2$ admit amenable actions then the free
product $G_1 \ast G_2$ (without amalgamation) also admits such an
action (J.L.\ Tu - private communication).  Using the equivalence with
the notion of exactness it follows that the same is true with
arbitrary amalgamation as well.  If $\gamma_i \in Aut(G_i)$ are
automorphisms then there is a natural free product automorphism
$\gamma_1 \ast \gamma_2 \in Aut(G_1 \ast G_2)$.

\begin{question}
Do we have that $ha(\gamma_1 \ast \gamma_2) = max(ha(\gamma_1),
ha(\gamma_2))$?
\end{question}

Note that if both $G_1, \ G_2$ are finite groups then $ha(\gamma_1
\ast \gamma_2) = 0$ since $\gamma_1 \ast \gamma_2$ is of finite order.

\section{Crystallographic Groups}

In this section we will compute the dual entropy of an arbitrary
automorphism of an extension of a finitely generated free abelian
group by a finite group.  However, we first require some elementary
preliminary results regarding F$\o$lner sets in finitely generated
abelian groups.

Fix $p \in {\Bbb N}$, let $v_1, \ldots, v_p$ be linearly independant
vectors in ${\Bbb R}^p$ and let $\chi = \{ v_1, \ldots, v_p \}$.
For each $t \in (0, \infty)$, let $\Gamma_{\chi}(t) \subset {\Bbb
R}^p$ be the parallelpiped $$\{ \sum_{i = 1}^{p} s_i v_i : s_i \in
{\Bbb R}, |s_i| \leq t, 1 \leq i \leq p \}.$$ For any two sets $H, K$
we let $H \triangle K$ denote the symmetric difference.

\begin{lemma}
For each $\delta > 0$ there exists a constant $C = C(\delta)$ with the
following property: For any set of linearly independant vectors $\chi
= \{ v_1, \ldots, v_p \}$ such that $$\Gamma_{\chi}(1) \supset \{ h
\in {\Bbb R}^p: h = (h_1, \ldots, h_p), |h_i| \leq 1, 1 \leq i \leq p
\},$$ (where we are using the canonical basis on the right hand side)
and any $x \in \Gamma_{\chi}(1) \cap {\Bbb Z}^p$ we have $$\frac{|(x +
\Gamma_{\chi}(C) \cap {\Bbb Z}^p) \triangle (\Gamma_{\chi}(C) \cap
{\Bbb Z}^p)|}{|\Gamma_{\chi}(C) \cap {\Bbb Z}^p|} < \delta.$$
\end{lemma}

{\noindent\bf Proof.}  Let $\chi = \{ v_1, \ldots, v_p \}$ be a basis
of ${\Bbb R}^p$ with the property stated above and let $C > 3$ be
chosen so that $(\frac{C - 2}{C + 1})^p > 1 - \delta/2$.

Letting $Q(1) = \Gamma_{\chi}(1) \cap {\Bbb Z}^p$ and $x \in Q(1)$ be
arbitrary we have $$\frac{|(x + Q(C)) \cap Q(C)|}{|Q(C)|} \geq
\frac{|Q(C-1)|}{|Q(C)|},$$ since $Q(C-1) \subset (x + Q(C)) \cap
Q(C)$.  However, our hypothesis on the set
$\chi$ implies the inequalities $|Q(C)| \leq vol(\Gamma_{\chi}(C +
1))$ and $vol(\Gamma_{\chi}(C - 2)) \leq |Q(C-1)|$, where $vol(\cdot)$
denotes the usual volume in ${\Bbb R}^p$.  Hence $$\frac{|(x + Q(C))
\cap Q(C)|}{|Q(C)|} \geq \bigg(\frac{C - 2}{C + 1}\bigg)^p.$$ But this
implies $$\frac{|(x + \Gamma_{\chi}(C) \cap {\Bbb Z}^p) \triangle
(\Gamma_{\chi}(C) \cap {\Bbb Z}^p)|}{|\Gamma_{\chi}(C) \cap {\Bbb
Z}^p|} \leq 2(1 - \bigg(\frac{C - 2}{C + 1}\bigg)^p) < \delta.  \ \ \
\ \ \ \Box$$

We will need to recall [Vo, Lem.\ 5.2]. 

\begin{lemma}
Let $\gamma \in GL(p, {\Bbb R})$ with eigenvalues $\lambda_1, \ldots,
\lambda_p$ and define $\mu_i = max\{ 1, |\lambda_i| \}$ for $1 \leq i
\leq p$.  There exists a basis $\{ v_1, \ldots, v_p \}$ of
${\Bbb R}^p$ with the property that if $\varepsilon > 0$ and a finite
subset $\sigma \subset {\Bbb R}^p$ are given then there exists $n_0
\in {\Bbb N}$ such that for all $n > n_0$ we have $$\{ \gamma^j(h) : h
\in \sigma, 1 \leq j \leq n \} \subset \{ \sum_{i = 1}^{p} s_i (1 +
\varepsilon)^n \mu_i^n v_i : s_i \in {\Bbb R}, |s_i| \leq 1, 1 \leq i
\leq p \}.$$
\end{lemma}

Assume $G = {\Bbb Z}^p \oplus F$ where $F$ is a finite abelian group
and let $\sigma \in Aut(G)$.  Then $\sigma(F) = F$ and hence there is
an induced automorphism $\dot{\sigma} \in Aut(G/F) = GL(p, {\Bbb Z})$.
Let $\lambda_1, \ldots, \lambda_p$ be the eigenvalues of
$\dot{\sigma}$ and $\mu_i = max\{ 1, |\lambda_i| \}$ for $1 \leq i
\leq p$.

\begin{lemma}
Let $G$, $\sigma \in Aut(G)$ and $\mu_1, \ldots, \mu_p$ be as above.
For each finite set $\omega \subset G$, $\delta > 0$ and $\varepsilon
> 0$ there exists a constant $K = K(\omega, \delta, \dot{\sigma})$ such that
$$ra(\omega + \cdots + \sigma^n(\omega), 2\delta) \leq |F| (Kn)^p(1 +
\varepsilon)^{np}(\mu_1 \mu_2 \cdots \mu_p)^n.$$
\end{lemma}

{\noindent\bf Proof.}  Letting $\pi : G \to G/F \cong {\Bbb Z}^p$ be
the quotient map, we see from Lemma 5.11 that $$ra(\omega + \cdots +
\sigma^n(\omega), 2\delta) \leq |F| ra(\pi(\omega) + \cdots +
\dot{\sigma}^n(\pi(\omega)), \delta).$$ 

Now choose a basis $\chi = \{ v_1, \ldots, v_p \}$ of ${\Bbb R}^p$ as
in Lemma 6.2.  Then we can find a constant $K_1 = K_1 (\omega, \chi)$
such that for all $n \in {\Bbb N}$ we have the inclusion $$\bigcup_{j
= 0}^{n} \dot{\sigma}^j (\pi(\omega)) \subset \{ \sum_{i = 1}^{p} s_i
(1 + \varepsilon)^n \mu_i^n v_i : s_i \in {\Bbb R}, |s_i| \leq K_1, 1
\leq i \leq p \}.$$ (The constant $K_1$ simply insures that the first
few iterates of $\pi(\omega)$ are also contained in the right hand
side.)  Hence we have $$\pi(\omega) + \cdots +
\dot{\sigma}^n(\pi(\omega)) \subset \{ \sum_{i = 1}^{p} s_i (1 +
\varepsilon)^n \mu_i^n v_i : s_i \in {\Bbb R}, |s_i| \leq nK_1, 1 \leq
i \leq p \},$$ for all $n \in {\Bbb N}$. 

Now let $\chi_n = \{ nK_1(1 + \varepsilon)^n \mu_i^n v_i \}_{1 \leq i
\leq p}$ and $\Gamma_{\chi_n} (\cdot)$ be as in Lemma 6.1.  We may
assume without loss of generality that each $\chi_n$ satisfies the
technical requirements of that lemma and hence find a constant $C = C
(\delta)$ such that for all $n \in {\Bbb N}$ and any $x \in
\Gamma_{\chi_n}(1) \cap {\Bbb Z}^p$, $$\frac{|(x + \Gamma_{\chi_n}(C)
\cap {\Bbb Z}^p) \triangle (\Gamma_{\chi_n}(C) \cap {\Bbb
Z}^p)|}{|\Gamma_{\chi_n}(C) \cap {\Bbb Z}^p|} < \delta.$$ Now defining
$T_n \in l^1({\Bbb Z}^p)$ to be the characteristic function over
$\Gamma_{\chi_n}(C) \cap {\Bbb Z}^p$ and $S_n = |\Gamma_{\chi_n}(C)
\cap {\Bbb Z}^p|^{-1} T_n$ one checks that $\| x.S_n - (S_n) \| <
\delta$ for all $x \in \Gamma_{\chi_n}(1) \cap {\Bbb Z}^p$.  In
particular  this shows $ra(\pi(\omega) + \cdots +
\dot{\sigma}^n(\pi(\omega)), \delta) \leq |supp(S_n)|$.  However, we
also have that $|supp(S_n)| = |\Gamma_{\chi_n}(C) \cap {\Bbb Z}^p|
\leq vol(\Gamma_{\chi_n}(C + 1)) = \tilde{C} (CK_1 n)^p(1 +
\varepsilon)^{np}(\mu_1 \mu_2 \cdots \mu_p)^n,$ where $\tilde{C}$ is
some constant depending only on $\chi$.  Letting $K =
\max(1,\tilde{C})CK_1$ concludes the proof of the lemma.  $\Box$

Assume now that $G$ contains a finitely generated normal abelian
subgroup $A$ of finite index.  Then $A \cong {\Bbb Z}^p \oplus F$ for
some finite abelian group $F$ and, as above, any automorphism of
$\sigma \in Aut(A)$ will induce an automorphism $\dot{\sigma} \in
GL(p, {\Bbb Z})$.  Also as above we let $\mu_j = \max(1,
|\lambda_j|)$, $1 \leq j \leq p$, for the eigenvalues $\lambda _j$ of
$\dot{\sigma}$.

\begin{theorem}
If $G$ is as above and $\gamma \in Aut(G)$ is such that $\gamma$
restricts to an automorphism $\sigma \in Aut(A)$ then $ha(\gamma) =
ha(\sigma) = ha(\dot{\sigma}) = \log(\mu_1 \mu_2 \cdots \mu_p)$.
\end{theorem}

{\noindent\bf Proof.}  First note that by duality we have
$ha(\dot{\sigma}) = \log(\mu_1 \mu_2 \cdots \mu_p)$ (cf.\ Theorem 4.2).
So by monotonicity and Corollary 5.12 it suffices to show $\log(\mu_1
\mu_2 \cdots \mu_p) \geq ha(\gamma).$ Let $\dot{\gamma}$ denote the
induced automorphism of $G/A$.  Since $ha(\gamma^k) = kha(\gamma)$ for
all positive integers (Proposition 5.3), it suffices to prove the
corresponding inequality for some power of $\gamma$.  Hence, replacing
$\gamma$ with a suitable power, we further assume that $\dot{\gamma} =
id \in Aut(G/A)$. 

Now choose a unital splitting $\phi : G/A \to G$ and let $\theta : G/A
\times G/A \to A$ be the associated cocycle map.  For each element $x
\in G/A$ let $a_x \in A$ be the unique element such that
$\gamma(\phi(x)) = \phi(x)a_x$ (since $\dot{\gamma}$ is the identity).
Note that for all integers $j$ we have $\gamma^j(\phi(x)) =
\phi(x)a_x\cdots\gamma^{j - 2}(a_x)\gamma^{j - 1}(a_x)$. Hence if
$\omega = \{ (\phi(x_1)a_1)^{-1}, \ldots , (\phi(x_m)a_m)^{-1} \}
\subset G$ is an arbitrary finite set and we define $\chi_1 = \{ a_1,
\ldots, a_m \} \cup \{ a_x : x \in G/A \} \subset A$ then Lemma 5.11
implies that $ra( \omega \cup \ldots \cup \gamma^n (\omega),
2\delta)$ is bounded above by $$|G/A| ra(\big(\bigcup_{x \in G/A}
\phi(x)^{-1} (\chi_1 + \gamma(\chi_1) + \cdots + \gamma^n (\chi_1) )
\phi(x)\big) + \big( \bigcup_{s,t \in G/A} \theta(s,t) \big), \delta),$$
for all $n \in {\Bbb N}$.

Using the fact that convolution is abelian when the underlying group
is abelian, one checks that $ra(X + Y,\delta) \leq ra(X,
\delta/2)ra(Y, \delta/2)$ for all finite subsets $X, Y \subset A$. So
letting $C = |G/A|ra(\cup_{s,t \in G/A} \theta(s,t) \big), \delta/2)$
we have the inequality $$ra( \omega \cup \ldots \cup \gamma^n
(\omega), 2\delta) \leq C ra(\bigcup_{x \in G/A} \phi(x)^{-1}
(\chi_1 + \gamma(\chi_1) + \cdots + \gamma^n (\chi_1) ) \phi(x),
\delta/2),$$ for all $n \in {\Bbb N}$.  Finally, define $\chi_2 =
\cup_{x \in G/A} \phi(x)^{-1} \chi_1 \phi(x) \subset A$ and we have
$$ra( \omega \cup \ldots \cup \gamma^n (\omega), 2\delta) \leq
Cra(\chi_2 + \ldots + \gamma^n (\chi_2), \delta/2),$$ since $\cup_{x
\in G/A} \phi(x)^{-1}(\gamma^j(X))\phi(x) = \gamma^j (\cup_{x \in G/A}
\phi(x)^{-1} X \phi(x) )$ for all finite subsets $X \subset A$ and $j
\in {\Bbb N}$ (this uses that $A$ is abelian). 

The previous lemma now states that for every $\varepsilon > 0$ we can
find some constant $K$ such that $$ra(\chi_2 + \ldots + \gamma^n
(\chi_2), \delta/2) \leq K n^p (1 + \varepsilon)^{np} (\mu_1 \mu_2
\cdots \mu_p)^n,$$ for all $n \in {\Bbb N}$ and hence
\begin{multline*} 
\begin{aligned} 
ha(\gamma, \omega, 2\delta)
    &= \limsup_{n \to \infty} \frac{\log(ra( 
       \omega \cup \ldots \cup \gamma^n (\omega), 2\delta))}{n + 1} \\[2mm]   
    &\leq \limsup_{n \to \infty} \frac{\log(CKn^p (1 + 
       \varepsilon)^{np} (\mu_1 \mu_2 \cdots \mu_p)^n)}{n + 1} \\[2mm] 
    &= \log(1 + \varepsilon)^p + \log(\mu_1 \mu_2 \cdots \mu_p).
\end{aligned}
\end{multline*}
But since $\varepsilon$ is arbitrary, this proves the theorem.  $\Box$

\begin{remark}
{\em A proof similar to that given above shows that if $1 \to A \to G \to H
\to 1$ is exact with $A$ finitely generated and abelian then for each
$\gamma \in Aut(G)$ such that $\gamma (A) = A$ and $\dot{\gamma} \in
Aut(H)$ is the identity we have $ha(\gamma) = ha(\gamma|_A)$.  In
particular, if $H$ is also abelian then the dual entropy of all inner
automorphisms of $G$ is determined by the restrictions to
$A$. (Compare with Corollary 5.15.)  }
\end{remark}

 If $H \subset G$ is a subgroup then we let $S_H = \{ g \in G : gh =
hg$ for all $h \in H \}$ denote the stabilizer of $H$ and $Z(H)$
denote the center.

\begin{lemma} 
Assume $A \subset G$ is a normal torsion free abelian subgroup of
finite index.  Then every automorphism of $G$ leaves $Z(S_A)$
invariant (i.e. $\gamma(Z(S_A)) = Z(S_A)$ for all $\gamma \in
Aut(G)$).
\end{lemma} 

{\noindent\bf Proof.}  Set $d = |G/A|$ and consider the set $G^d = \{
g^d : g \in G \}$.  Note that $G^d \subset A$ and if $\gamma \in
Aut(G)$ then $\gamma(G^d) = G^d$.  Thus if $g \in \gamma(S_A)$ then
$g$ commutes with $a^d$ for all $a \in A$ since $a^d \in G^d =
\gamma(G^d) \subset \gamma(A)$.  Thus $a^d = g a^d g^{-1} =
(gag^{-1})^d$.  But, since $A$ is abelian, this implies $1_g =
(a^{-1}gag^{-1})^d$ which implies $ga = ag$ since $A$ is torsion
free.  Thus $\gamma(S_A) \subset S_A$ and the same argument applied
to $\gamma^{-1}$ implies $\gamma(S_A) = S_A$.  But if $\gamma$ leaves
$S_A$ invariant then it also leaves $Z(S_A)$ invariant.  $\Box$ 

We are finally in a position to compute $ha(\cdot)$ for all
automorphisms of crystallographic groups.  First, assume that $1 \to
{\Bbb Z}^p \to G \to F \to 1$ is exact, with $F$ a finite group.  Let
$A = Z(S_{{\Bbb Z}^p})$.  Then $A$ is a finitely generated normal
abelian subgroup of $G$ of finite index.  Writing $A \cong {\Bbb Z}^q
\oplus L$ for some $q \in {\Bbb N}$ and finite abelian group $L$ we
can define a homomorphism $\rho : Aut(G) \to GL(q,{\Bbb Z})$ in the
following manner. If $\gamma \in Aut(G)$ then the previous lemma says
that $\gamma$ leaves $A$ invariant and hence $\gamma|_A$ leaves $L$
invariant and hence defines an element $\rho(\gamma) \in Aut(A/L) =
GL(q,{\Bbb Z})$. As a  consequence of Theorem 6.4 we get:

\begin{corollary}
Let $G$ and $\rho : Aut(G) \to GL(q,{\Bbb Z})$ be as above.  Then for
every $\gamma \in Aut(G)$, $$ha(\gamma) = ha(\rho(\gamma)) =
\log(\mu_1 \mu_2 \cdots \mu_q),$$ where $\mu_j = max(1, |\lambda_j|)$,
$1 \leq j \leq q$, for each eigenvalue $\lambda_j$ of $\rho(\gamma)$.
\end{corollary}

{\noindent\bf Acknowledgements.}  The first named author would like to
thank his thesis advisor, Marius Dadarlat, for a number of important
discussions related to entropy in (amenable) groups.  Indeed, a number
of the ideas presented here can be traced back to those conversations. 

The present collaboration was made possible by the authors'
participation in the program Probabilit\'es Libres et Espaces
D'Op\'erateurs, held at Institut Henri Poincar\'e.  We wish to thank
the institute for it's support and hospitality and the organizers for
putting together such a stimulating semester.

\end{document}